\documentclass[11pt]{amsart}    

\usepackage{graphicx}
\usepackage{color}
\usepackage{enumerate}
\usepackage{amsfonts,amsmath,latexsym,amssymb,amsthm}
\usepackage{hyperref}
\usepackage{geometry}
\usepackage{txfonts}
\usepackage{color}
\definecolor{Myred}{cmyk}{0.0,1.0,1.0,0.00}

\newtheorem{claim}{Claim}[section]
\newtheorem{theorem}[claim]{Theorem}
\newtheorem{lemma}[claim]{Lemma}

\begin{document}

\title[Magnetic Dirichlet Laplacian in curved waveguides]
{Magnetic Dirichlet Laplacian in curved waveguides}

\author{Diana Barseghyan} 
\address{Department of Mathematics, University of Ostrava, 30. dubna 22, 70103 Ostrava, Czech Republic}
\email{diana.schneiderova@osu.cz}
\author{Swanhild Bernstein}
\address{Institute of Applied Analysis, TU Bergakademie Freiberg\\
Akademiestrasse 6, Freiberg 09599, Germany}
\email{Swanhild.Bernstein@math.tu-freiberg.de}
\author{Baruch Schneider}
\address{Department of Mathematics, University of Ostrava, 30. dubna 22, 70103 Ostrava, Czech Republic}
\email{baruch.schneider@osu.cz}
\author{Martha Lina Zimmermann}
\address{Institute of Applied Analysis, TU Bergakademie Freiberg\\
Akademiestrasse 6, Freiberg 09599, Germany}
\email{Martha-Lina.Zimmermann@math.tu-freiberg.de}


\maketitle

\begin{abstract}
For a two-dimensional curved waveguide, it is well known that the
spectrum of the Dirichlet Laplacian is unstable. Any perturbation of the straight strip produces eigenvalues below the essential spectrum. In this paper, a magnetic field is added. We explicitly prove that the spectrum of the magnetic Laplacian is stable under small but non-local deformations of the waveguide.
\end{abstract}
\bigskip


\bigskip

\section{Introduction} \label{s: intro}
\setcounter{equation}{0}

It has long been known that an appropriate bending of a two
dimensional quantum waveguide induces the existence of bound states,
\cite{ES89}, \cite{GJ92} and \cite{DE95}. From a mathematical point of
mathematical point of view, this means that the Dirichlet Laplacian on a smooth asymptotically straight planar waveguide has at least one isolated eigenvalue below the threshold of the essential spectrum.
Similar results have been obtained for a locally deformed
waveguide, which corresponds to adding a small ``bump'' to the straight
waveguide, see \cite{BGRS97} and \cite{BEGK01}. As a result
at least one isolated eigenvalue appears below the essential spectrum
for any nonzero curvature, satisfying certain regularity
properties.

As is well-known, that a magnetic field, even a local one, can significantly affect the behaviour of waveguide systems, in particular the existence of a geometrically induced discrete spectrum. While a particle confined in a fixed-profile tube with a hard-wall boundary can exist in localized states whenever the tube is bent or locally deformed (and asymptotically straight), cf. \cite{EK15} for a comprehensive review of quantum waveguide theory, the presence of a local magnetic field can destroy such a discrete spectrum.

The analogous effect of bound state existence resulting from the geometry of the interaction support has been observed is a class of singular Schr\"{o}dinger operators, usually dubbed leaky quantum wires, with attractive contact interaction supported by a curve \cite{EK15}. It was established in \cite{BE21} that the presence of a local magnetic field can again destroy such a discrete spectrum. Another similar result has been established in \cite{BBS24}, where the authors consider the magnetic Schr\"{o}dinger operator with a non-negative potential supported over a strip which is a local deformation of a straight one, and the magnetic field is assumed to be nonzero and local. 
The object of our interest in this paper will be the magnetic Dirichlet Laplacian, i.e., the self-adjoint operator  
\begin{equation}\label{1}H_\Omega(A) = (i\nabla + A)^2\end{equation}  
associated with the quadratic form 
$$\|(i\nabla + A)u\|_{L^2(\Omega)},\quad u \in \mathcal{H}_0^1(\Omega),$$
where the real-value function $A$ is a vector potential.

Our main results are presented in Section 3. Here we prove the stability of the essential spectrum of the operator (\ref{1}) with respect to the situation when the strip is a smooth asymptotically straight planar waveguide and the magnetic field is non-zero and vanishing at infinity. We also give a sufficient condition for the absence of the discrete spectrum.
Let us mention that the similar results have been obtained in \cite{EK05}, but the authors in that work deal only with compactly supported perturbations of the straight tube. In the present work it has been shown that the vanishing and non-zero magnetic field can destroy the discrete spectrum even if the perturbation of the tube is not local.

\section{Curved waveguides} \label{s: twodim}
\setcounter{equation}{0}

The domain $\Omega$ studied in our paper is assumed to be a curved planar strip of a width $d$. Its points are described by the curvilinear coordinates $s\in\mathbb{R}, u\in(0, d)$ as follows:
\begin{eqnarray}\nonumber x=x(s, u)= a(s) -u \dot{b}(s),\\\label{coordinates}
y =y(s, u)= b(s) + u \dot{a}(s),
\end{eqnarray}
where dot marks the derivative with respect to $s\,$ and $a,b$ are smooth functions characterizing the reference curve 
$\Gamma = \{(a(s), b(s)): s\in\mathbb{R}\}$. We assume
\begin{equation}\label{ab}
\dot{a}(s)^2 + \dot{b}(s)^2 = 1,
\end{equation}
so $s$ is the arc length of $\Gamma$, and $u\in(0, d)$ is the distance of the point $(x, y)$ from $\Gamma$. It is useful to introduce the signed curvature $\gamma(s)$ of $\Gamma$,
\begin{equation}
\label{curv.}
\gamma(s)=\dot{b}(s)\ddot{a}(s)- \dot{a}(s)\ddot{b}(s).
\end{equation}

The characterization of the region (\ref{coordinates}) by the curvilinear coordinate only makes sense if the latter can be uniquely defined, which imposes two different restrictions. First, the transverse size must not be too large, the inequality $d |\gamma(s)| < 1$ must hold at every point of the curve. In addition, the region $\Omega$ does not intersect itself. For obvious reasons, $\Omega$ is called a curved quantum waveguide.

Let us create the following identity that we will use later. Given (\ref{curv.}) and (\ref{ab}) we have
\begin{eqnarray}\nonumber\label{gamma}\gamma(s)^2= \left(\dot{b}(s)\ddot{a}(s)- \dot{a}(s)\ddot{b}(s)\right)^2\\\nonumber
=\ddot{a}(s)^2 \dot{b}(s)^2+ \ddot{b}(s)^2\dot{a}(s)^2- 2\ddot{a}(s)\dot{b}(s)\ddot{b}(s)\dot{a}(s)\\\nonumber
=\ddot{a}(s)^2(1-\dot{a}(s)^2) + \ddot{b}(s)^2\dot{a}(s)^2- 2\ddot{a}(s)\dot{b}(s)\ddot{b}(s)\dot{a}(s)\\\label{aux}=
\ddot{a}(s)^2- \ddot{a}(s)^2 \dot{a}(s)^2 + \ddot{b}(s)^2\dot{a}(s)^2- 2\ddot{a}(s)\dot{b}(s)\ddot{b}(s)\dot{a}(s).
\end{eqnarray}

With respect to identity (\ref{ab}), we have 
\begin{equation}\label{adot}
\dot{a}(s)\ddot{a}(s)+\dot{b}(s)\ddot{b}(s)=0.
\end{equation}
Therefore, (\ref{aux}) implies 
\begin{eqnarray}\nonumber\gamma(s)^2=
\ddot{a}(s)^2- \ddot{b}(s)^2 \dot{b}(s)^2 + \ddot{b}(s)^2\dot{a}(s)^2+ 2\ddot{b}(s)^2\dot{b}(s)^2= \\\label{id}
=\ddot{a}(s)^2+\ddot{b}(s)^2.
\end{eqnarray}

During our work will need the following conditions 
\begin{enumerate}[(a)]
\setlength{\itemsep}{3pt}
\item $$|\gamma(s)|,\,|\dot{\gamma}(s)|,\,|\ddot{\gamma}(s)|,\,|\dddot{\gamma}(s)|\le \frac{\beta}{1+s^2},\quad s\in\mathbb{R},$$
where $\beta>0$ is some constant to be described later.

\item Suppose that the magnetic vector potential $A(x, y) = (a_1(x, y), a_2(x, y))$ is such that for $j = 1, 2$ the functions $a_j,\, \frac{\partial a_j}{\partial x},\, \frac{\partial a_j}{\partial y}$ vanish for large values of $x$.
\end{enumerate}

\section{Main results} \label{s: intro}
\setcounter{equation}{0}

We will establish the stability of the essential spectrum and the absence of the discrete spectrum below the threshold of the essential spectrum of the operator (\ref{1}) defined on an asymptotically straight curved waveguide under the assumptions (a) and (b). Let us start with the absence of the discrete spectrum due to a non-zero magnetic field.

\subsection{Absence of the discrete spectrum}
\begin{theorem}\label{absence}
Let $B \in C^1(\mathbb{R}^2)$ be a real-valued magnetic field which is non-trivial in $\Omega$. Assume that assumptions (a) and (b) hold. There exists a positive number $\beta_0$ such that for $\beta\in (0, \beta_0)$ the discrete spectrum of the operator (\ref{1}) below $\frac{\pi^2}{d^2}$ is empty. 
\end{theorem}

\begin{proof}
Let us denote by $q$ the quadratic form of the operator (\ref{1}). Then
\begin{equation}\label{qf}
q_\Omega^A(\varphi)=\int_\Omega\left(\left|i\frac{\partial\varphi}{\partial x}+ a_1\varphi\right|^2+\left|i\frac{\partial\varphi}{\partial y}+ a_2\varphi\right|^2\right)\,d x\,d y,
\end{equation}
where $A=(a_1, a_2)$.

Define the unitary operator
\begin{equation} \label{U1}
U : L^2(\Omega) \to L^2(\Omega_0), \quad \Omega_0= \mathbb{R}\times (0, d), 
\end{equation}
which for any $\varphi\in L^2(\Omega)$ acts as
\begin{equation} \label{U2}
\psi(s, u):=(U \varphi)(s, u) = \sqrt{1 + u \gamma(s)}\, 
\varphi(a(s) - u \dot{b}(s), b(s) + u \dot{a}(s)).
\end{equation}

Moreover, if $\varphi|_{\partial\Omega}=0$, then the same is true for $\psi$: $\psi|_{\partial\Omega_0}=0$.

Let us now follow the calculations of \cite{EK05}. We can check that the Jacobian
\begin{equation}\label{part.}
\frac{\partial(x, y)}{\partial(s, u)}= 1+u \gamma(s)
\end{equation}
and
\begin{eqnarray}\nonumber
\frac{\partial\varphi}{\partial x} = (1 +u \gamma)^{-1} \left(\dot{a}\frac{\partial }{\partial s} - (\dot{b} +
u \ddot{a}) \frac{\partial}{\partial u}\right)\left(\frac{\psi}{\sqrt{1 + u
\gamma}} \right) \\\label{part.1} \frac{\partial\varphi}{\partial y} = (1 + u \gamma)^{-1}
\left( \dot{b} \frac{\partial}{\partial s} - (\dot{a}- u \ddot{b}) \right)\left(\frac{\psi}{\sqrt{1 + u
\gamma}} \right).
\end{eqnarray}

Then, using the notation 
\begin{equation} \label{tildeA}
\tilde A(s, u) = \left(\tilde{a_1}(s, u), \tilde{a_2}(s, u)\right) = \left(a_1\left(a(s) - u \dot{b}(s), b(s) +u \dot{a}(s)\right), a_2\left(a(s) - u \dot{b}(s), b(s) +u \dot{a}(s)\right)\right)
\end{equation}
we rewrite the quadratic form $q_\Omega^A(\varphi),\, \varphi\in\mathcal{H}_0^1(\Omega)$ as follows
\begin{eqnarray} \nonumber
q_\Omega^A(\varphi)= \int_{\Omega_0} \left( \left|\left(\frac{i}{1 + u \gamma}\left( \dot{a}
\frac{\partial}{\partial s} - (\dot{b} + u \ddot{a})\frac{\partial}{\partial u} \right) + \tilde a_1 \right) \left( \frac{\psi} {\sqrt{1 + u
\gamma}} \right) \right|^2 \right. \nonumber  \\\label{q0}
 \left. \quad + \left| \left(\frac{i}{1+u \gamma} \left( \dot{b} \frac{\partial}{\partial s}
+ (\dot{a} - u \ddot{b}) \frac{\partial}{\partial u}\right) +
\tilde a_2 \right) \left( \frac{\psi}{\sqrt{1 + u
\gamma}} \right)\right|^2 \right) (1 + u \gamma) \,d s\,d u.
\end{eqnarray}

Using (\ref{adot}) and (\ref{id}) expression (\ref{q0}) performs 
\begin{gather}\nonumber
q_\Omega^A(\varphi)= \int_{\Omega_0} \biggl(\frac{1}{(1+ u\gamma)^2}\left|\frac
{\partial\psi}{\partial s}\right|^2 + \frac{i (\dot{a}\tilde
{a_1} + \dot{b} \tilde{a_2})}{1 + u \gamma} \left(\frac{\partial\psi}{\partial s} \overline
{\psi} - \psi \frac{\partial \overline{\psi}}{\partial s}\right) + \left|\frac{\partial\psi}{\partial u}\right|^2\\ \nonumber
+ \frac{i \left( - (\dot{b}+ u \ddot{a}) \tilde{a_1} + (\dot{a} - u \ddot{b})
\tilde{a_2}\right)}{1 + u \gamma} \left(\frac{\partial\psi}{\partial u} \overline{\psi}
-\psi  {\frac{\partial\overline{\psi}}{\partial u}} \right) \\\nonumber
- \frac{u \dot{\gamma} }{2(1 + u \gamma)^3}\left(\psi
{\frac{\partial\overline{\psi}}{\partial s}} + \frac{\partial\psi}{\partial s} \overline{\psi}\right) - \frac{\gamma}{2(1 + u \gamma)} \left(\psi{\frac{\partial \overline{\psi}}{\partial u}} +
\frac{\partial\psi}{\partial u} \overline{\psi}\right)\\  \label{gamma}
+ \left( \frac{u^2  \left(\dot{\gamma} \right)^2}
{4(1 + u \gamma)^4} + \frac{ \gamma^2}{4(1 + u
\gamma)^2} + \tilde{a_1}^2 + \tilde{a_2}^2\right) |\psi|^2\biggr)
\,d s\,d u\,.
\end{gather}

We write the right-hand side of (\ref{gamma}) as a perturbation of the form $\tilde{q}_\Omega^A$ as follows
\begin{equation}
\label{pert.}
q_\Omega^A(\varphi)= \tilde{q}_\Omega^A- I(\psi)\,,
\end{equation}
where
\begin{eqnarray}\nonumber
\tilde{q}_\Omega^A(\psi)= \int_{\Omega_0} \left|i \frac{\partial\psi}{\partial s} + (\dot{a}
\tilde{a_1} + \dot{b} \tilde{a_2}) \psi\right|^2 + \left|i \frac{\partial\psi}{\partial u} +
 (-\dot{b} \tilde{a_1} + \dot{a} \tilde{a_2}) \psi\right|^2 \,d s\, d u\,,\\
\nonumber
I(\psi)
= \int_{\Omega_0} \biggl( \frac{2 u \gamma + u^2
\gamma^2}{1 + u \gamma} \left|\frac{\partial\psi}{\partial s}\right|^2 + i u \gamma (\dot{a} \tilde{a_1} + \dot{b} \tilde{a_2}) \left(\frac{\partial\psi}{\partial s} \overline\psi
- \psi{\frac{\partial \overline{\psi}}{\partial s}}\right)\\\nonumber+ i u \left(-  \gamma \dot{b} \tilde{a_1} +
\gamma \dot{a} \tilde{a_2} + \frac{\ddot{a} \tilde{a_1} - \ddot{b} \tilde{a_2}} {1 +
u \gamma}\right) \left( \frac{\partial\psi}{\partial u} \overline \psi -
\psi {\frac{\partial \overline{\psi}}{\partial u}} \right)  \\\nonumber
+ \frac{u  \dot{\gamma}}{2(1 + u \gamma)^3}
\left(\psi{\frac{\partial \overline{\psi}}{\partial s}}+ \frac{\partial\psi}{\partial s} \overline \psi
\right) + \frac{\gamma}{2(1 + u \gamma)} \left(\psi
{\frac{\partial \overline{\psi}}{\partial u}} + \frac{\partial\psi}{\partial u} \overline \psi \right)
- \left( \frac{u^2 \left( \dot{\gamma} \right)^2}
{4(1 + u \gamma)^4} + \frac{ \gamma^2} {4(1 + u
\gamma)^2} \right) |\psi|^2 \biggr) \,d s \,d u\,,
\end{eqnarray}

Now let us estimate $I(\psi)$. We will use the notation $\|\cdot\|_\infty= \|\cdot\|_{L^\infty(\mathbb{R})}$.
Since $u\in (0, d)$, then given (\ref{ab}) and (\ref{id}) we get  
\begin{eqnarray}\nonumber
| I(\psi)|\le \frac{2d+ d^2\|\gamma_\infty\|}{1-d\|\gamma\|_\infty} \int_{\Omega_0}|\gamma|\left|\frac{\partial\psi}{\partial s}\right|^2\,d s\,d u
\\\nonumber +2d\sqrt{\|\tilde{a_1}\|_\infty^2+\|\tilde{a_2}\|_\infty^2} \int_{\Omega_0}|\gamma||\psi|\left|\frac{\partial\psi}{\partial s}\right|\,d s\,d u+ 2d\sqrt{\|\tilde{a_1}\|_\infty^2+\|\tilde{a_2}\|_\infty^2}\int_{\Omega_0}|\gamma||\psi|\left|\frac{\partial\psi}{\partial u}\right|\,d s\,d u
\\\nonumber+ \frac{2d \sqrt{\|\tilde{a_1}\|_\infty^2+\|\tilde{a_2}\|_\infty^2} }{1-d\| \gamma\|_\infty}\int_{\Omega_0}|\gamma||\psi|\left|\frac{\partial\psi}{\partial u}\right|\,d s\,d u+\frac{d}{(1-d \|\gamma\|_\infty)^3}\int_{\Omega_0}|\dot{\gamma}||\psi|\left|\frac{\partial\psi}{\partial s}\right|\,d s\,d u\\\nonumber+ \frac{1}{1-d\| \gamma\|_\infty}\int_{\Omega_0}|\gamma||\psi|\left|\frac{\partial\psi}{\partial u}\right|\,d s\,d u+\frac{d^2}{4(1-d \|\gamma\|_\infty)^4}\int_{\Omega_0}|\dot{\gamma}|^2|\psi|^2\,d s\,d u\\\nonumber+
\frac{1}{4(1-d \|\gamma\|_\infty)^2}\int_{\Omega_0}\gamma^2|\psi|^2\,d s\,d u
\le\frac{2d+ d^2\|\gamma_\infty\|}{1-d\|\gamma\|_\infty} \int_{\Omega_0}|\gamma|\left|\frac{\partial\psi}{\partial s}\right|^2\,d s\,
d u\\\nonumber +d\sqrt{\|\tilde{a_1}\|_\infty^2+\|\tilde{a_2}\|_\infty^2} \int_{\Omega_0}|\gamma|\left|\frac{\partial\psi}{\partial s}\right|^2\,d s\,d u+ d\sqrt{\|\tilde{a_1}\|_\infty^2+\|\tilde{a_2}\|_\infty^2} \int_{\Omega_0}|\gamma||\psi|^2\,d s\,d u\\\nonumber
+ d\sqrt{\|\tilde{a_1}\|_\infty^2+\|\tilde{a_2}\|_\infty^2}\int_{\Omega_0}|\gamma|\left|\frac{\partial\psi}{\partial u}\right|^2\,d s\,d u
+ d\sqrt{\|\tilde{a_1}\|_\infty^2+\|\tilde{a_2}\|_\infty^2}\int_{\Omega_0}|\gamma||\psi|^2\,d s\,d u
\\\nonumber+ \frac{d \sqrt{\|\tilde{a_1}\|_\infty^2+\|\tilde{a_2}\|_\infty^2} }{1-d\| \gamma\|_\infty}\int_{\Omega_0}|\gamma|\left|\frac{\partial\psi}{\partial u}\right|^2\,d s\,d u+ \frac{d \sqrt{\|\tilde{a_1}\|_\infty^2+\|\tilde{a_2}\|_\infty^2} }{1-d\| \gamma\|_\infty}\int_{\Omega_0}|\gamma||\psi|^2\,d s\,d u\\\nonumber+\frac{d}{2(1-d \|\gamma\|_\infty)^3}\int_{\Omega_0}|\dot{\gamma}|\left|\frac{\partial\psi}{\partial s}\right|^2\,d s\,d u +\frac{d}{2(1-d \|\gamma\|_\infty)^3}\int_{\Omega_0}|\dot{\gamma}||\psi|^2\,d s\,d u\\\nonumber+ \frac{1}{2(1-d\| \gamma\|_\infty)}\int_{\Omega_0}|\gamma|\left|\frac{\partial\psi}{\partial u}\right|^2\,d s\,d u+
\frac{1}{2(1-d\| \gamma\|_\infty)}\int_{\Omega_0}|\gamma||\psi|^2\,d s\,d u
\\\nonumber+\frac{d^2}{4(1-d \|\gamma\|_\infty)^4}\int_{\Omega_0}|\dot{\gamma}|^2|\psi|^2\,d s\,d u+
\frac{1}{4(1-d \|\gamma\|_\infty)^2}\int_{\Omega_0}\gamma^2|\psi|^2\,d s\,d u\\\nonumber
\le\left(\frac{d (2+ d\|\gamma\|_\infty)}{1-d\|\gamma\|_\infty}+ d\sqrt{\|\tilde{a_1}\|_\infty^2+\|\tilde{a_2}\|_\infty^2}+\frac{d}{2(1-d \|\gamma\|_\infty)^3}\right) \int_{\Omega_0}\left(|\gamma|+ |\dot{\gamma}|\right)\left|\frac{\partial\psi}{\partial s}\right|^2\,d s\,d u 
\\\nonumber
+\left(d\sqrt{\|\tilde{a_1}\|_\infty^2+\|\tilde{a_2}\|_\infty^2}+ \frac{2d \sqrt{\|\tilde{a_1}\|_\infty^2+\|\tilde{a_2}\|_\infty^2} +1}{2(1-d\| \gamma\|_\infty)}\right) \int_{\Omega_0}|\gamma|\left|\frac{\partial\psi}{\partial u}\right|^2\,d s\,d u
\\\nonumber +\left(2d \sqrt{\|\tilde{a_1}\|_\infty^2+\|\tilde{a_2}\|_\infty^2}+ \frac{2d \sqrt{\|\tilde{a_1}\|_\infty^2+\|\tilde{a_2}\|_\infty^2}+1}{2(1-d\|\gamma\|_\infty)}+ \frac{d}{2(1-d\|\gamma\|_\infty)^3}+\frac{d^2}{4(1-d \|\gamma\|_\infty)^4}+ \frac{1}{4(1-d \|\gamma\|_\infty)^2}\right)\\\nonumber\int_{\Omega_0}\left(|\gamma|+ \gamma^2+ |\dot{\gamma}|+ \dot{\gamma}^2\right) |\psi|^2\,d s\,d u\,.
\end{eqnarray}

Let us denote
\begin{equation}\label{alpha1}
\alpha_1:=\mathrm{max}\left\{\frac{d (2+ d\|\gamma\|_\infty)}{1-d\|\gamma\|_\infty}+ d\sqrt{\|\tilde{a_1}\|_\infty^2+\|\tilde{a_2}\|_\infty^2}+\frac{d}{2(1-d \|\gamma\|_\infty)^3}, \,d\sqrt{\|\tilde{a_1}\|_\infty^2+\|\tilde{a_2}\|_\infty^2}+ \frac{2d \sqrt{\|\tilde{a_1}\|_\infty^2+\|\tilde{a_2}\|_\infty^2}+1 }{2(1-d\| \gamma\|_\infty)}\right\}\,,
\end{equation}
\begin{equation}\label{alpha2}
\alpha_2:=2d \sqrt{\|\tilde{a_1}\|_\infty^2+\|\tilde{a_2}\|_\infty^2}+ \frac{2d \sqrt{\|\tilde{a_1}\|_\infty^2+\|\tilde{a_2}\|_\infty^2}+1}{2(1-d\|\gamma\|_\infty)}+ \frac{d}{2(1-d\|\gamma\|_\infty)^3}+\frac{d^2}{4(1-d \|\gamma\|_\infty)^4}+ \frac{1}{4(1-d \|\gamma\|_\infty)^2}\,.
\end{equation}

Hence the above inequality implies 
$$ |I(\psi)|\le \alpha_1\int_{\Omega_0}\left(|\gamma|+ |\dot{\gamma}|\right) |\nabla \psi|^2\,d s\,d u+
 \alpha_2\int_{\Omega_0}\left(|\gamma|+ \gamma^2+ |\dot{\gamma}|+ \dot{\gamma}^2\right) |\psi|^2\,d s\,d u\,.
$$

In view of (\ref{pert.}) and the above estimate, we arrive
\begin{equation}
\label{est.}
q_\Omega^A(\varphi)\ge \tilde{q}_\Omega^A- \alpha_1\int_{\Omega_0}\left(|\gamma|+|\dot{\gamma}|\right)|\nabla \psi|^2\,d s\,d u-
 \alpha_2\int_{\Omega_0}\left(|\gamma|+ \gamma^2+|\dot{\gamma}|+ \dot{\gamma}^2\right) |\psi|^2\,d s\,d u\,.
\end{equation}

One can easily check that for any magnetic potential $\hat{A}= (\hat{a_1}, \hat{a_2})$ the following pointwise inequality holds
$$
\left|\nabla \psi\right|^2\le 2\left|i\nabla \psi+ \hat{A}\psi\right|^2 + 2\left(\hat{a_1}^2+ \hat{a_2}^2\right)
|\psi|^2\,, \quad \psi\in \mathcal{H}^1(\Omega_0)\,.
$$

In view of this fact and (\ref{est.}) we get
\begin{equation}
\label{est.1}
q_\Omega^A(\varphi)\ge \int_{\Omega_0}\left(1-\rho_1^\gamma\right)\left|i \nabla \psi+ \tilde{A} \psi\right|^2\,d s\,d u- \int_{\Omega_0}\rho_2^\gamma |\psi|^2\,d s\,d u\,,
\end{equation}
where 
\begin{eqnarray}\nonumber
\rho_1^\gamma= 2\alpha_1 \left(|\gamma|+ |\dot{\gamma}|\right)\,,
\\\label{rho}
\rho_2^\gamma= 2\alpha_1\left(\tilde{a_1}^2+ \tilde{a_2}^2\right) \left(|\gamma|+|\dot{\gamma}|\right)+
\alpha_2 \left(|\gamma|+ \gamma^2+ |\dot{\gamma}|+ \dot{\gamma}^2\right)\,.
\end{eqnarray}

For the further proof we need the following Hardy-type inequality.
\begin{lemma}\label{lemma1}
For any function $g\in \mathcal{H}_0^1(\Omega_0)$, the following estimate holds    
 \begin{eqnarray*}
\int_{\Omega_0}\left(f^2 \left|i \nabla g+ \tilde{A} g\right|^2- \frac{\pi^2}{d^2} f^2|g|^2\right)\,d s\,d u\ge C_{\tilde{A}}\int_{\Omega_0}\frac{f^2}{1+s^2} |g|^2\,d s\,d u+\int_{\Omega_0}f'' |g|^2\,d s\,d u\,,
\end{eqnarray*}
where $C_{\tilde{A}}>0$ is a constant, $\tilde{A}=(\tilde{a_1}, \tilde{a_2})$ is the magnetic potential and $f= f(s):\mathbb{R}\to (0, \infty)$ is a smooth function. 
\end{lemma}

Let us return to the proof of the lower bound (\ref{est.1}). Using the lemma \ref{lemma1} with $f(s)=\left(1-\rho_1^\gamma(s)\right)^{1/2}$ we get
\begin{eqnarray}\nonumber
q_\Omega^A(\varphi)- \frac{\pi^2}{d^2} \int_{\Omega_0}\left(1-\rho_1^\gamma(s)\right)|\psi|^2\,d s\,d u\ge C_{\tilde{A}}\int_{\Omega_0}\frac{\left(1-\rho_1^\gamma\right)}{1+s^2} |\psi|^2\,d s\,d u
- \frac{1}{4}\int_{\Omega_0}\frac{\left(2(1-\rho_1^\gamma) \ddot{\rho_1^\gamma}+ (\dot{\rho_1^\gamma})^2\right)}{\left(1-\rho_1^\gamma\right)^{3/2}}
|\psi|^2\,d s\,d u\\\label{est.2}- \int_{\Omega_0}\rho_2^\gamma |\psi|^2\,d s\,d u\,.
\end{eqnarray}

Let us choose a positive number $\beta_0$ in assumption (a)  such that for any $\beta\in (0, \beta_0)$ we have
$$C_{\tilde{A}} \frac{\rho_1^\gamma}{1+s^2}+ \frac{1}{4}\frac{\left(2(1-\rho_1^\gamma) \ddot{\rho_1^\gamma}+(\dot{\rho_1^\gamma})^2\right)}{\left(1-\rho_1^\gamma\right)^{3/2}}
+ \frac{\pi^2 \rho_1^\gamma}{d^2}+ \rho_2^\gamma \le  \frac{C_{\tilde{A}}}{1+s^2}\,,\quad s\in \mathbb{R}\,.
$$

This establishes that 
$$
q_\Omega^A(\varphi)- \frac{\pi^2}{d^2} \int_{\Omega_0}|\psi|^2\,d s\,d u= q_\Omega^A(\varphi)- \frac{\pi^2}{d^2} \int_{\Omega_0}|\varphi|^2\,d s\,d u\ge 0\,,
$$
which proves that the discrete spectrum below $\frac{\pi^2}{d^2}$ of the operator (\ref{1}) is empty.
\end{proof}

\subsection{Stability of the essential spectrum}
\begin{theorem}\label{stability}
Suppose that the magnetic vector potential $A = (a_1, a_2)$ is such that for $j = 1, 2$ the functions $a_j,\, \frac{\partial a_j}{\partial x},\, \frac{\partial a_j}{\partial y}$ are vanishing for large values of $x$. Then the essential spectrum of the operator (\ref{1}) coincides with the half-line $\left[\frac{\pi^2}{d^2}, \infty\right)$.
\end{theorem}

\begin{proof}

To prove that any non-negative number $\mu\ge \frac{\pi^2}{d^2}$ belongs to the essential
spectrum of $H_\Omega(A)$, we will use Weyl's criterion
\cite[Thm.~VII.12]{RS81}: we have to find a sequence
$\{\varphi_n\}_{n=1}^\infty\subset D(H_\Omega(A))$ of unit vectors,
$\|\varphi_n\|=1$, which converges weakly to zero and
$$
\|H_\Omega(A)\varphi_n-\mu\varphi_n\|\to 0 \qquad\text{as}\quad n\to\infty
$$
holds. 

First, let us rewrite operator (\ref{1}) in the following form
$$
H_\Omega(A)= -\Delta+ 2i \left(a_1 \frac{\partial}{\partial x}+ a_2 \frac{\partial}{\partial y} \right)+ i\left(\frac{\partial a_1}{\partial x}+ \frac{\partial a_2}{\partial y}\right)+ \left(a_1^2+ a_2^2\right).
$$
Then
\begin{eqnarray}\nonumber
\int_\Omega\left|H_\Omega(A)\varphi_n- \mu \varphi_n\right|^2\,d x\,d y\\\nonumber=
\int_\Omega\biggl|-\Delta \varphi_n+ 2i \left(a_1 \frac{\partial\varphi_n}{\partial x}+ a_2 \frac{\partial\varphi_n}{\partial y} \right)+ i\left(\frac{\partial a_1}{\partial x}+ \frac{\partial a_2}{\partial y}\right)\varphi_n+ \left(a_1^2+ a_2^2\right)\varphi_n- \mu\varphi_n\biggr|^2\,d x\,d y
\end{eqnarray}
and 
\begin{eqnarray}\nonumber
\int_\Omega\left|H_\Omega(A)\varphi_n- \mu \varphi_n\right|^2\,d x\,d y\\\nonumber\le
4\int_\Omega\left|-\Delta \varphi_n- \mu\varphi_n\right|^2\,d x\,d y+ 16\int_\Omega\left|a_1 \frac{\partial\varphi_n}{\partial x}+ a_2 \frac{\partial\varphi_n}{\partial y}\right|^2\,d x\,d y\\\label{Weyl}+ 4\int_\Omega\left|\frac{\partial a_1}{\partial x}+ \frac{\partial a_2}{\partial y}\right|^2|\varphi_n|^2\,d x\,d y+ 4\int_\Omega \left(a_1^2+ a_2^2\right)^2|\varphi_n|^2\,d x\,d y\,.
\end{eqnarray}

For each $\mu= \frac{\pi^2}{d^2}+ k^2,\,k\in\mathbb{Z}$, we will separately estimate each integral in (\ref{Weyl}).
Let us mention that, using the transformation (\ref{U1}), it was established in \cite{ES89} that the Dirichlet Laplacian $-\Delta$ on $\Omega$ can be rewritten as the operator $H_0$ on $L^2(\Omega_0)$, which acts as
$$
(H_0\psi)(s,u) = -\frac{\partial}{\partial s}\left(\frac{1}{(1+u\gamma(s))^2}\frac{\partial\psi}{\partial s}\right)-\frac{\partial^2\psi}{\partial u^2}+W(s, u) \psi\,,
$$
$$
W(s,u):= -\frac{\gamma^2(s)}{4(1+u\gamma(s))^2}+\frac{u\ddot\gamma(s)}{2(1+u\gamma(s))^3}-\frac{5}{4}
\frac{u^2\dot\gamma^2(s)}{(1+u\gamma(s))^4}
$$
and Dirichlet boundary conditions are imposed at $u=0,\, d$. 

Given the unitary equivalence, it is sufficient to estimate the first integral in (\ref{Weyl}) for the operator $H_0$.
We will construct the Weyl sequence $\varphi_n\in C_0^\infty(\Omega)$ in such a way that the corresponding mapping (\ref{U1}) sequence $\psi_n\in C_0^\infty(\Omega_0)$ is the following normalized sequence  
\begin{equation}\label{psi}
\psi_n(s, u)= \sqrt{\frac{2}{d n}}f\left(\frac{s}{n}\right) e^{i k s} \sin\left(\frac{\pi u}{d}\right),\quad n\in \mathbb{N},
\end{equation}
where $f\in C_0^\infty(\mathbb{R})$ is a smooth function with support in the interval $(1, 2)$ with $L^2$ norm equal to one.

We have 
\begin{eqnarray}\nonumber
\int_{\Omega_0}\left|H_0\psi_n- \mu\psi_n\right|^2\,d x\,d y\\\nonumber= \int_{\Omega_0}\left| -\frac{\partial}{\partial s}\left(\frac{1}{(1+u\gamma(s))^2}\frac{\partial \psi_n}{\partial s}\right)-\frac{\partial^2\psi_n}{\partial u^2}+W(s, u) \psi_n- \mu\psi_n \right|
^2\,d s\,d u\\\nonumber
= \frac{2}{d n}\int_{\Omega_0}\biggl|\frac{2 u \dot{\gamma}}{(1+u \gamma)^3} \left(\frac{1}{n}\dot{f}\left(\frac{s}{n}\right)+ i k f\left(\frac{s}{n}\right)\right)- \frac{1}{(1+u \gamma)^2}\biggl(\frac{1}{n^2} \ddot{f}\left(\frac{s}{n}\right)+ \frac{2 i k}{n} \dot{f}\left(\frac{s}{n}\right)- k^2 f\left(\frac{s}{n}\right)\biggr) \\\nonumber+ \frac{\pi^2}{d^2} f\left(\frac{s}{n}\right) 
+W(s, u) f\left(\frac{s}{n}\right)- \mu f\left(\frac{s}{n}\right) \biggr|
^2 \sin^2\left(\frac{\pi u}{d}\right)\,d s\,d u\\\nonumber
=\frac{2}{d n}\int_{\Omega_0}\biggl|\frac{2 u \dot{\gamma}}{(1+u \gamma)^3} \left(\frac{1}{n}\dot{f}\left(\frac{s}{n}\right)+ i k f\left(\frac{s}{n}\right)\right)- \frac{1}{(1+u \gamma)^2}\biggl(\frac{1}{n^2} \ddot{f}\left(\frac{s}{n}\right)+ \frac{2 i k}{n} \dot{f}\left(\frac{s}{n}\right)\biggr) \\\nonumber
-\frac{k^2(2u \gamma+u^2\gamma^2)}{(1+u \gamma)^2} f\left(\frac{s}{n}\right) 
+W(s, u) f\left(\frac{s}{n}\right)\biggr|^2 \sin^2\left(\frac{\pi u}{d}\right)\,d s\,d u\,.\end{eqnarray}

Because of assumption (a), it is easy to check that the right-hand side of the above estimate is $\mathcal{O}\left(\frac{1}{n^2}\right)$. Therefore  
\begin{equation}\label{H0}
\int_{\Omega_0}\left|H_0\psi_n- \mu\psi_n\right|^2\,d x\,d y= \mathcal{O}\left(\frac{1}{n^2}\right)\,.
\end{equation}

Let us pass to the second term in the right-hand side of (\ref{Weyl}) and estimate as the first $\int_\Omega a_1^2 \left|\frac{\partial\varphi_n}{\partial x}\right|^2\,d x\,d y$. We pass to the curvilinear coordinates and use (\ref{part.}), (\ref{part.1}) and the notations from (\ref{tildeA}). With the notation 
$$
\omega_n=\left\{(x(s, u), y(s, u))_{s\in(n, 2n),\, u\in(0, d)}\right)\}$$
one gets
\begin{eqnarray}\nonumber&&
\int_\Omega a_1^2 \left|\frac{\partial\varphi_n}{\partial x}\right|^2\,d x\,d y
=\int_{\Omega_0}\tilde{a_1}^2\left|(1+ u\gamma)^{-1}\left(\dot{a} \frac{\partial}{\partial s}- (\dot{b}+u \ddot{a})\frac{\partial}{\partial u}\right)\left(\frac{\psi_n}{\sqrt{1+ u\gamma}}\right)\right|^2 (1+ u \gamma)\,d s\,d u\\\nonumber&&\le
\frac{\|\tilde{a_1}\|_{L^\infty(\omega_n)}^2}{(1-d \|\gamma\|_\infty)^2}\int_{\Omega_0} \left|\frac{\dot{a}}{\sqrt{1+ u\gamma}} \frac{\partial \psi_n}{\partial s}- \frac{u \dot{a} \dot{\gamma}}{2(1+u \gamma)^{3/2}} \psi_n-\frac{(\dot{b}+ u\ddot{a})}{\sqrt{1+u \gamma}}\frac{\partial \psi_n}{\partial u}+ \frac{(\dot{b}+ u\ddot{a}) \gamma}{2(1+ u\gamma)^{3/2}}\psi_n\right|^2 (1+u \gamma)\,d s\,d u\\\nonumber&&\le
\frac{4 \|\tilde{a_1}\|_{L^\infty(\omega_n)}^2}{(1-d \|\gamma\|_\infty)^2} \int_{\Omega_0}
\biggl(\frac{\dot{a}^2}{(1- d\|\gamma\|_\infty)}\left|\frac{\partial \psi_n}{\partial s}\right|^2 + \frac{d^2 \dot{a}^2 \dot{\gamma}^2}{4(1- d \|\gamma\|_\infty)^3} |\psi_n|^2\\\nonumber&&+ \frac{2(\dot{b}^2+ d^2\ddot{a}^2)}{1-d \|\gamma\|_\infty}\left|\frac{\partial \psi_n}{\partial u}\right|^2+ \frac{(\dot{b}^2+ d^2\ddot{a}^2)}{2(1- d\|\gamma\|_\infty)^3}|\psi_n|^2\biggr) (1+ d\|\gamma\|_\infty)\,d s\,d u\,.
\end{eqnarray}
 
In view of (\ref{ab}) and (\ref{id}) and the construction of $\psi_n$ 
\begin{equation}\label{last}
\int_\Omega a_1^2 \left|\frac{\partial\varphi_n}{\partial x}\right|^2\,d x\,d y=\mathcal{O}\left(\|\tilde{a_1}\|_{L^\infty(\omega_n)}^2\right) \int_{\Omega_0}\left(\left|\frac{\partial \psi_n}{\partial s}\right|^2+ \left|\frac{\partial \psi_n}{\partial u}\right|^2+ |\psi_n|^2\right)\,d s\,d u= \mathcal{O}\left(\|\tilde{a_1}\|_{L^\infty(\omega_n)}^2\right)\,.
\end{equation}

In the same way one can estimate the integrals
$$\int_\Omega a_2^2 \left|\frac{\partial\varphi_n}{\partial y}\right|^2\,d x\,d y=\mathcal{O}\left(\|\tilde{a_2}\|_{L^\infty(\omega_n)}^2\right)\,.
$$

For the remaining third and fourth integrals, we have
\begin{eqnarray*}\int_\Omega\left|\frac{\partial a_1}{\partial x}+ \frac{\partial a_2}{\partial y}\right|^2|\varphi_n|^2\,d x\,d y=
\mathcal{O}\left(\left\|\frac{\partial a_1}{\partial x}\right\|^2_{L^\infty(\omega_n)}+ \left\|\frac{\partial a_2}{\partial y}\right\|^2_{L^\infty(\omega_n)}\right)\,,\\
\int_\Omega \left(a_1^2+ a_2^2\right)^2|\varphi_n|^2\,d x\,d y= \mathcal{O}\left(\|\tilde{a_1}\|^4_{L^\infty(\omega_n)}+\|\tilde{a_2}\|_{L^\infty(\omega_n)}^4\right)\,.
\end{eqnarray*}

Combining the last four expressions, (\ref{H0}) and the inequality (\ref{Weyl}) one establishes our claim.
\end{proof}

\section{Proof of Lemma \ref{lemma1}}

Let us denote $$h(s, u):= f(s) g(s, u).$$ Then
\begin{eqnarray}\nonumber
\int_{\Omega_0} \left|i \nabla h+ \tilde{A} h\right|^2\,d s\,d u\\\nonumber= \int_{\Omega_0}\left|i f \frac{\partial g}{\partial s} +
 \tilde{a_1} f g +i f' g \right|^2\,d s\,d u
+ \int_{\Omega_0} \left|i f \frac{\partial g}{\partial u} + \tilde{a_2} f g\right|^2\,d s\,d u
\\\nonumber=
\int_{\Omega_0}f^2 \left|i \frac{\partial g}{\partial s} + \tilde{a_1} g\right|^2\,d s\,d u
-i \int_{\Omega_0}\left(i \frac{\partial g}{\partial s} +\tilde{a_1} g \right)f f' \overline{g}\,d s\,d u\\\nonumber
+i \int_{\Omega_0}\left(-i{\frac{\partial  \overline{g}}{\partial s}} +\tilde{a_1} \overline{g} \right)f f' g\,d s\,d u+
\int_{\Omega_0} f'^2 |g^2|\,d s\,d u +\int_{\Omega_0}f^2 \left|i  \frac{\partial g}{\partial u} + \tilde{a_2}  g\right|^2\,d s\,d u\\\nonumber
=\int_{\Omega_0}f^2 \left|i \frac{\partial g}{\partial s} + \tilde{a_1} g\right|^2\,d s\,d u
+\int_{\Omega_0}\frac{\partial g}{\partial s} f f' \overline{g}\,d s\,d u\\\nonumber
+\int_{\Omega_0}\frac{\partial \overline{g}}{\partial s}  f f' g\,d s\,d u
+\int_{\Omega_0} f'^2 |g^2|\,d s\,d u +\int_{\Omega_0}f^2 \left|i  \frac{\partial g}{\partial u} + \tilde{a_2}  g\right|^2\,d s\,d u\\\nonumber
=\int_{\Omega_0}f^2 \left|i \nabla g+ \tilde{A} g\right|^2\, ds\,d u+\int_{\Omega_0}\frac{\partial g}{\partial s} f f' \overline{g}\,d s\,d u\\\label{lemma}+\int_{\Omega_0}\frac{\partial \overline{g}}{\partial s}  f f' g\,d s\,d u
+\int_{\Omega_0} f'^2 |g^2|\,d s\,d u\,.
\end{eqnarray}

One can check that 
$$
\int_{\Omega_0}\frac{\partial g}{\partial s} f f' \overline{g}\,d s\,d u+ \int_{\Omega_0}\frac{\partial \overline{g}}{\partial s}  f f' g\,d s\,d u\\
=- \int_{\Omega_0} (f'^2+ f'') |g|^2\,d s\,d u\,.
$$

Hence the right-hand side of (\ref{lemma}) becomes
\begin{equation}\label{last}
\int_{\Omega_0} \left|i \nabla h+ \tilde{A} h\right|^2\,d s\,d u=\int_{\Omega_0}f^2 \left|i \nabla g+ \tilde{A} g\right|^2\, ds\,d u- \int_{\Omega_0} f'' |g|^2\,d s\,d u \,.
\end{equation}

Now we are going to use the following Hardy-type inequality, see e.g. \cite{EK05}
\begin{theorem}
Let $B \in C^1(\mathbb{R}^2)$ be a real-valued magnetic field such that $B\not\equiv 0$ in $D= \mathbb{R}\times(0, \pi)$. Then 
$$
c_B \int_D\frac{|u|^2}{1+s^2}\, d s\,d u \le \int_D\left(\left|i \nabla u +A u\right|^2- |u|^2\right)\, d s\,d u\,, 
$$
holds for all $u \in\mathcal{H}_0^1(D)$, where $A$ is a magnetic vector potential associated with $B$ and  $c_B$ is a positive constant.
\end{theorem}

The above theorem can be transformed as follows
 \begin{theorem}\label{theorem}
Let $B \in C^1(\mathbb{R}^2)$ be a real-valued magnetic field such that $B\not\equiv 0$ in $\Omega_0= \mathbb{R}\times(0, d)$. Then 
$$
\tilde{c}_B \int_{\Omega_0}\frac{|u|^2}{1+s^2}\, d s\,d u \le \int_{\Omega_0}\left(\left|i \nabla u +A u\right|^2- \frac{\pi^2}{d^2} |u|^2\right)\, d s\,d u\,, 
$$
holds for all $u \in\mathcal{H}_0^1(\Omega_0)$, where $A$ is a magnetic vector potential associated with $B$ and the constant in the left-hand side is represented via the constant from the theorem above as follows   
$$\tilde{c}_B=c_{\frac{d^2}{\pi^2} B\left(\frac{s d}{\pi}, \frac{u d}{\pi}\right)}.$$ 
\end{theorem}

Employing above theorem at the left-hand side of (\ref{last}) one obtains 
\begin{eqnarray}\nonumber
\tilde{c}_{\tilde{B}} \int_{\Omega_0}\frac{|h|^2}{1+s^2}\, d s\,d u+ \frac{\pi^2}{d^2} \int_{\Omega_0}|h|^2\,d s\,d u\le
\int_{\Omega_0}f^2 \left|i \nabla g+ \tilde{A} g\right|^2\, ds\,d u- \int_{\Omega_0}  f'' |g|^2\,d s\,d u \,,
\end{eqnarray}
where $\tilde{B}$ is the magnetic field corresponding to the magnetic potential $\tilde{A}$.
This establishes the statement of the Lemma \ref{lemma1} with the constant $C_{\tilde{A}}= \tilde{c}_{\tilde{B}}$.

\end{document}